\documentclass[12pt]{article}
\input{epsf}

\usepackage{amsmath,amssymb,amsfonts,graphicx}
\makeatletter \@addtoreset{equation}{section} \makeatother
\newtheorem{thm}{Theorem}[section]
\newtheorem{lemma}{Lemma}[section]
\newtheorem{prop}{Proposition}[section]
\newtheorem{cor}{Corollary}[section]

\newtheorem{rem}{Remark}[section]

\begin{document}
\begin{center}

{\bf \large On the Distribution of Products of Spherical Classes
in Classical Symmetric Spaces of Rank One}

\vspace*{1cm} {\bf Jafar Shaffaf}\footnote{Institute for Studies
in Theoretical Physics and Mathematics (IPM) and Sharif University
of Technology, Tehran, Iran.  Email : shaffaf@ipm.ir}
\vspace{0.3cm}
\end{center}

\vspace*{1cm}

\begin{abstract}
\noindent  The distribution of products of random matrices chosen
from fixed spherical classes is determined for classical rank 1
symmetric spaces.  It is observed that $n\to\infty$ limit behaves
approximately as in the abelian case. A theorem on the rate of
convergence to the Haar measure in the case of $SU(n)$ is also
established.
\end{abstract}

\begin{center}
\end{center}

\section{Introduction}
A basic problem in random matrix theory is the determination of
the distribution of products of matrices chosen randomly from
specified classes.  For example the problem of the support of
products of random unitary matrices chosen from fixed conjugacy
classes, or sums of Hermitian matrices of given eigenvalues is
treated in [AW] and [F] (which is based on the work of Klyachko
[Kl]). In this paper we investigate the distribution of products
of random matrices chosen from spherical classes for classical
rank 1 symmetric spaces.  When matrices are chosen from conjugacy
or spherical classes from a simple group of $2\times 2$ matrices a
complete solution is given in [Sh1] and [Sh2].

The distribution of products of spherical classes may be
interpreted as the determination of the algebra structure for the
Hecke algebra generated by the singular measures concentrated on
orbits of $K\times K$ in $G$ for a symmetric pair $(G,K)$. It can
also be given an interpretation in geometric probability. Consider
for example $S^{2n-1}$ with the natural action of $SU(n)$ on it.
Fix a longitude $\Lambda\subset S^{2n-1}$, let $z_i\in \Lambda$,
$i=1,2$, and choose matrices $g_i$ randomly from $SU(n)$
(according to Haar measure) conditioned on the requirement that
$g_i.z_i$ lies in the same meridian as $z_i$. Then one may inquire
about the density function for $g_2g_1$.

For the case of classical symmetric spaces of rank 1 a complete
solution to this problem is given in this paper (Theorems
\ref{thm:shaffaf4}, \ref{thm:shaffaf5}, \ref{thm:shaffaf23} and
\ref{thm:shaffaf24}, \ref{thm:shaffaf22}). It is observed that as
the dimension of the symmetric space of a given class tends to
infinity, the distribution of the product measures converges
weakly to a singular continuous measure concentrated on a single
orbit of $K\times K$.  This may be rephrased as the $n\to\infty$
limit behaves approximately as in the abelian case. A similar
limiting behavior is also observed for higher rank symmetric
spaces but this more complex limit theorem will be treated later
in another publication.

In view of this limiting behavior, in section 7 the rate of
convergence to Haar measure on $SU(n)$ of products of the form
\begin{eqnarray*}
g_1h_1g_2h_2\ldots g_Nh_N,
\end{eqnarray*}
where $g_i$'s (resp. $h_i$'s) are chosen from a fixed spherical
class ${\mathcal O}_a$ (resp. ${\mathcal O}_b$), is determined.
In particular, it is shown that for $N \sim C\log n$ the product
measure $\lambda_a\star\lambda_b\star\ldots \star
\lambda_a\star\lambda_b$, (product of $N$ copies of
$\lambda_a\star\lambda_b$) where $\lambda_a$ is the (probability)
invariant measure on ${\mathcal O}_a$, tends to the Haar measure
of $SU(n)$ as $n\to\infty$ in $L^p$ for $1\le p\le 2$.

The methods used in this work are based on harmonic analysis on
symmetric spaces.  In section 2 relevant integration formulae for
symmetric spaces are stated in the appropriate form.  Sections 3,
4 , 5 and 6 give the explicit density functions for products of
spherical classes in classical symmetric spaces of rank 1 of both
compact and non-compact types.  In the final section essential use
is made of representation theory and harmonic analysis for
$SU(n)/S(U(1)\times U(n-1))$ to obtain the rate of convergence.
Similar results are also valid for rank 1 symmetric spaces of the
orthogonal groups but are not treated here.

The author wishes to thank Professor S. Shahshahani who gave him
the opportunity to pursue his interests, and especially Professor
Mehrdad Shahshahani for suggesting the problem and many
stimulating discussions.

\section{Integration Formulae}

In this section we recall some basic integration formulae related
to symmetric spaces.  A detailed treatment is given in [H1] and
[H2].

Let $G$ be a connected compact semi-simple Lie group ,
$\mathfrak{g}$ its Lie algebra and $K$ be a subgroup of $G$ with
Lie algebra $\mathfrak{k}$ such that $(G, K)$ is a symmetric pair
of compact type.  Let $ \mathfrak{g}= \mathfrak{k} +  \mathfrak{p}
$ the corresponding Cartan decomposition, $\mathfrak{a} \subset
\mathfrak{p}$ a maximal abelian subspace and $\Sigma$ the
corresponding set of restricted roots. Fix a Weyl chamber
$\mathfrak{a}^+\subset\mathfrak{a}$, let $\Sigma^+$ and
$\Phi=\{\alpha_1,\ldots,\alpha_l\}\subset \Sigma_+$ denote the
corresponding sets of simple and positive roots. The multiplicity
of a root $\alpha$ will be denoted by $m_\alpha$.  Let $M$ be the
centralizer of $\mathfrak{a}$ in $\mathfrak{k}$, the group $A =
\exp\mathfrak{a}$ is a closed, and therefore a compact subgroup of
$G$.

Let $dg$, $dk$, $dm$ and $da$ denote the Haar measures on the
compact groups $G$, $K$, $M$ and $A$, and $du$ and $db$ be
invariant measures on $U=G/K$, $B=K/M$, respectively. We consider
the surjective map
\begin{eqnarray*}
\Psi: (K/M)\times A\longrightarrow G/K,~~~\Psi(kM, a)= kaK.
\end{eqnarray*}
The Jacobian of $\Psi$ is
\begin{equation}
\label{eq:shaffaf2}
 \det ( d\Psi_{(kM, a)} ) = \prod_{
\alpha  \in  \Sigma^+}  | \sin \alpha (i \ H)|^   {m_ {\alpha}}
\end{equation}
where $a = \exp (H)$, $H \in \mathfrak{a}$ and $m _\alpha$ is the
multiplicity of the restricted root $\alpha$.  Let $\delta(a)$
denote the right hand sight of (\ref{eq:shaffaf2}).  $\Psi$ is
one-to-one and regular on an open dense set and we have :
\begin{eqnarray*}
 \Psi ^\star (du)= c  \delta(a) db da
 \end{eqnarray*}
where $c$ is a constant depending on the normalization of measures.

\begin{thm}
\label{thm:shaffaf1} Let $U=G/K$ be a Riemannian symmetric space
of compact type. Then with the above notation we have
\begin{eqnarray*}
\int_G f ( gK) dg = c \int_{K/M} \left(\int_A f(ka K) \delta(a) da
\right)  db
\end{eqnarray*}
for all $f\in C(G/K)$. Moreover if $G$ is simply connected and $K$
is connected, then
\begin{eqnarray*}
 \int_{G/K} f(gK) d
u = c \int _ {K/M} db \left( \int_Q f( k ( \exp H)K) \delta( \exp
H) d H \right )
\end{eqnarray*}
where
\begin{eqnarray*}
\delta(\exp(H)) = \prod_{ \alpha  \in  \Sigma^+}  | \sin \alpha (i
\ H)|^ {m_ {\alpha}},~~~c^{-1} = \int_{Q} \prod_ {\alpha \in
\Sigma^+} ~ (\sin ( \alpha (i H) )^ {m_\alpha} d H,
\end{eqnarray*}
the measures on $U=G/K$, $B=K/M$ and $A$ are normalized to be 1,
and $Q$ is the polyhedron
\begin{eqnarray*}
Q = \{H\in \mathfrak{a}: \frac{1}{i} \mu_j(H)> 0, ( 1\leq j \leq l),
\frac{1}{i} \mu(H)< \pi \}.
\end{eqnarray*}
\end{thm}

According to the classification theory the classical symmetric
spaces of rank 1 are

\begin{center}
\renewcommand{\arraystretch}{1.2}
\begin{tabular}{|c|c|c|c|}\hline
Noncompact&Compact&Dimension of $G/K$\\\hline $(SU(n,1),S(U(1)\times
U(n)))$&$(SU(n+1), S(U(1)\times U(n)))$&$2n$\\\hline
$(SO(n,1),SO(n))$&$(SO(n+1),SO(n))$&$n$\\\hline $(Sp(n,1),Sp(1)
\times Sp(n))$&$(Sp(n+1),Sp(1)\times Sp(n))$&$4n$\\\hline
\end{tabular}
\end{center}
There is also an exceptional symmetric pair of rank 1 associated to the exceptional Lie group $F_4$ which is not treated
in this paper.

\section{The Symmetric Pair $(SU(n+1), S(U(1) \times U(n)))$ }

In this section $G= SU(n+1) $ and $K = S(U(1) \times U(n))$. In
this case $G$ is simply connected, $K$ connected and $G/K$ is
irreducible. The Lie algebra of $G$ is $\mathfrak{g}= su(n+1)$,
the space of traceless skew hermitian matrices and the Lie algebra
of the subgroup $K$ is
 \begin{eqnarray*}
 \mathfrak{k}=
\left\{\left[\begin{array}{c|ccc} -\text{tr}(A)&&0
\\\hline  &&\\0&&A\\&&
\end{array}\right]
 ; \  \textrm{$A$ is an $n\times n$ skew hermitian
 matrix}\right\}.
\end{eqnarray*}
Let $\mathfrak{g} = \mathfrak{k} \oplus \mathfrak{p}$ be the
Cartan decomposition then $\mathfrak{p}$ is the subspace of
matrices of the form
\begin{eqnarray*}
\left[\begin{array}{ccccc}
0& & \xi_1 & \ldots & \xi_n  \\
\begin{array}{c}
\\-\overline{\xi_1} \\
\vdots \\
\\-\overline{\xi_n} \\
\end{array} &  & & 0_{n \times n} & \\

\end{array}\right]
\end{eqnarray*}
A maximal abelian subspace $ \mathfrak{a} \subset \mathfrak{p}$ is
$ \mathfrak{a} = \{t H ~;~~ t \in R \}$ where $H$ is the matrix

 \begin{eqnarray*} H = \left[\begin{array}{ccccc}
0& & \ 1 &  \ldots &  0 \\
\begin{array}{c}
\ -1 \\
\ 0 \\
\vdots \\
\ 0 \\
\end{array} &  & & 0_{n \times n} & \\

\end{array}\right]
\end{eqnarray*}
The centralizer $\mathfrak{m}$ of $\mathfrak{a}$ in $\mathfrak{k}$
is the subalgebra of $\mathfrak{k}$ consisting of matrices of the
form
$$ \left[\begin{array}{c|ccc} -\text{tr}(A)&&0
\\\hline  &&\\0&&A\\&&
\end{array}\right]$$
where $A= [a_{ij}]$ is an $n\times n$ skew hermitian matrix
satisfying

\begin{eqnarray*}
A v  = \ \textrm{tr}(A) \ v
\end{eqnarray*}
where $v=[1,1, \ldots ,1]^T$. There is a unitary matrix $P$ such
that
\begin{eqnarray*}
PAP^{-1} \ = \ \left[\begin{array}{c|ccc} \textrm{tr}(A)&&0
\\\hline  &&\\0&&B\\&&
\end{array}\right],~~~{\rm tr}(B)=0.
\end{eqnarray*}
Therefore we can identify the Lie algebra $\mathfrak{m}$ with $s(
u(1) \times u(1) )\times su(n-1) $.

\begin{lemma}
The Jacobian $\delta$ for polar coordinates of  the symmetric pair
$(G,K)$ is given by the formula

\begin{eqnarray*}
\delta(\exp(tH))=(\sin t)^{2(n-1)}~\sin 2t .
\end{eqnarray*}
\end{lemma}

\noindent {\bf Proof}-
By a straightforward calculation the positive restricted roots are

\begin{eqnarray*}
 \alpha_1(tH) = i t , \ \ \alpha_2(t H ) = 2 \alpha_1 (tH)= 2i
t
\end{eqnarray*}
with multiplicities $ m _{\alpha_1} = 2(n-1), m_{\alpha_2} = 1 $.
Therefore by theorem (\ref{thm:shaffaf1}) for the Jacobian
$\delta$ we have

\begin{eqnarray*}
\delta(\exp(t H)) = \prod_{ \alpha  \in \  \Sigma^+}  | \sin
\alpha (i t H)|^ {m_ {\alpha}} = (\sin ( t) ) ^{2(n-1)} \ \sin (2
t),
\end{eqnarray*} and the lemma followed. $\blacksquare$

\begin{lemma}
The stabilizer subgroup $K_{z_0}$ of $z_0 = [0:1:1: \ldots :1]$,
under left action of $K$ on $\mathbb{CP}(n)\simeq G/K$ is the
centralizer subgroup $M$.
\end{lemma}
\noindent {\bf Proof}- Let $p=[z, \xi ] \in \mathbb{CP}(n)$  and
note that $\mathbb{CP}(n)$ is the quotient space $ \mathbb{CP}(n)=
S^{2n+1}/S^1$ under diagonal action of $S^1$.  The orbit of a
generic point $z=[\zeta, \xi ] \in \mathbb{CP}(n)$ under $K = S(
U(1) \times U(n))$ is

\begin{equation}
\label{eq:shaffaf51} \ \left[\begin{array}{c|ccc}
e^{-i\theta}&0&\cdots&0
\\\hline 0&&\\\vdots&&B\\0&&
\end{array}\right]\left[\begin{array}{cccc}
\zeta
\\\hline\\  \xi\\\\
\end{array}\right]=\left[\begin{array}{cccc}
e^{-i\theta}\zeta
\\\hline\\  B\xi\\\\
\end{array}\right].
\end{equation}
It is now clear that the first component of the above vector is a
copy of $S^1$ and the second component is a $(2n-1)$-sphere with
equivalency under $S^1$ which is, by definition, a copy of
$\mathbb{CP}(n-1)$.  Hence we showed that the orbit of a generic
point $z \in \mathbb{CP}(n)$ under the action of the subgroup $K$
is equivalent to $S^1 \times \mathbb{CP}(n-1)$.

Fix a generic point $[\zeta,\xi]\in \mathbb{CP}(n)$ and consider
the corresponding embedding of $S^1 \times
\mathbb{CP}(n-1)\hookrightarrow \mathbb{CP}(n)$ given by
(\ref{eq:shaffaf51}. Let $z \in \mathbb{CP}(n)$ we know $ O_z = K
\ .\ z = S^1 \times \mathbb{CP}(n-1)$ and let $K_z$ be the
stabilizer subgroup of the point $z$ \ i.e. $ K_z =\{ k \in K \ |
\ kz = z \}$.

We want to show that this stabilizer subgroup $K_z$ for $z=z_0$ as
in the theorem is isomorphic to the centralizer subgroup $M$.
Recall that the Lie algebra $\mathfrak{m}$ of the Lie group $M$ is
characterized by
$$ \mathfrak{m} = \left\{ B = \left[\begin{array}{c|ccc}
-\text{tr}(A)&&0
\\\hline  &&\\0&&A\\&&
\end{array}\right] \in \mathfrak{k} \ \Big {|} \ A^{\ast} = -A \ , \ B \ v =
\ \text{tr}(A) \ v \right\}
$$
where $v = [0, 1, 1, \ldots , 1]$.  So we obtain
$$ e^B \ v = \det( e^A) \ v. $$ Therefore $v$ is an eigenvector for $M$ and conversely every element in $K$ which has $v$ as
an eigenvector is necessarily an element of the centralizer
subgroup $M$.  Hence we have $K_{z_0}=M$. $\blacksquare$
\begin{rem}
\label{rem:shaffaf1} {\rm It is clear from the above lemma that the
centralizer subgroup $M$ is a connected subgroup and therefore we
have $M= S(U(1)\times U(1) \times U(n-1))$.}
\end{rem}

\begin{lemma}
We have $K/M \simeq S^1 \times \mathbb{CP}(n-1)$.
\end{lemma}
\noindent {\bf Proof}- According to the previous Lemma the orbit of the point $ z_0 = [0:1:1:
\ldots :1] \in \mathbb{CP}(n)$ under the action of the group $K$
can be identified with $S^1 \times \mathbb{CP}(n-1)$ and the stabilizer of this point
is the centralizer subgroup $M$, so by the isomorphism theorem we have
$$O_{z_0} = K/K_{z_0} =  K/M.$$
Hence $ K/M = S^1 \times \mathbb{CP}(n-1)$. $\blacksquare$ \\
Thus we are allowed to apply the polar coordinate on the pair $(K,
M)$ and we bring it as the following lemma
\begin{lemma}
The Jacobian of polar coordinates for the pair $(K,M)$ is
\begin{eqnarray*}
\delta_0(\exp(t H)) =  (\sin ( t) ) ^{2(n-2)} \ \sin
(2 t)
\end{eqnarray*}
\end{lemma}
\noindent {\bf Proof}- Since the pair $(K, M)$ is symmetric in
this case we have the Cartan decomposition $ \mathfrak{k} =
\mathfrak{m} \oplus \mathfrak{p}_0$ with $$ [\mathfrak{m},
\mathfrak{m}]\subseteq \mathfrak{m}, \ \ [\mathfrak{m},
\mathfrak{p}_0]\subseteq \mathfrak{p}_0 , \ \ [\mathfrak{p}_0,
\mathfrak{p}_0]\subseteq \mathfrak{m}. $$ Now $\mathfrak{m} =
\mathfrak{s}(\mathfrak{u}(1) \times \mathfrak{u}(1)) \times
\mathfrak{su}(n-1)$ which are the matrices of the form
 $$ \left[\begin{array}{c|c|ccc}
- \text{tr}(A)&0&0&\cdots&0
\\\hline 0 & \text{tr}(A) & 0 & \  \cdots \ \ 0\\\hline
0& 0 &&&\\\vdots&\vdots&& B &\\ 0& 0&&&\\
\end{array}\right],  \ \ B \in su(n-1)
$$and the vector space $\mathfrak{p}_0$ is matrices in the Lie algebra
$\mathfrak{k}$ of the form
$$ \left[\begin{array}{c|c|ccc} 0&0&0&\cdots&0
\\\hline 0&0&\xi_1&\cdots&\xi_{n-1}\\\hline
0&-\overline{\xi}_1&&&\\\vdots&\vdots&&0_{n-2}&\\0&-\overline{\xi}_{n-1}&&&\\
\end{array}\right]
$$  It is clear that the maximal
subspace is spanned by
$$ H_0 = \left[\begin{array}{c|c|ccc} 0&0&0&\cdots&0
\\\hline 0&0& 1 & \cdots & 0 \\\hline
0& -1   &&&\\\vdots&\vdots&&0_{n-2}&\\ 0 & 0 &&&\\
\end{array}\right]\in \mathfrak{p}_0.
$$  The positive restricted roots are

\begin{eqnarray*}
 \beta_1(tH_0) = i t , \ \ \beta_2(t H_0 ) = 2 \beta_1 = 2i
t
\end{eqnarray*}
with multiplicities $ m _{\beta_1} = 2(n-2), m_{\beta_2} = 1 $.
Therefore by theorem (\ref{thm:shaffaf1}) we have

\begin{eqnarray*}
\delta_0(\exp(t H_0)) = \prod_{ \beta  \in \  \Sigma^+}  | \sin
\alpha (i t H)|^ {m_ {\alpha}} = (\sin ( t) ) ^{2(n-2)} \ \sin (2
t).
\end{eqnarray*} The lemma is proved. $\blacksquare$

Now we state the main theorem of this section.

\begin{thm}
\label{thm:shaffaf4} Let $\lambda_a$ and $\lambda_b$ be two
(singular) spherical measures concentrated on the $K$-spherical
classes ${\cal O}_a$ and ${\cal O}_b$ in the group $G= SU(n+1)$
respectively. Then $\lambda_a\star\lambda_b$ is absolutely
continuous relative to the Haar measure on $SU(n+1)$. It is a
spherical measure and for a continuous spherical function $f$ on
$SU(n+1)$ we have

\begin{eqnarray*}
\lambda_a\star\lambda_b(f) = (c\mathrm{vol}(M)\mathrm{vol}(K))^2
\delta (t_1) \delta (t_2) \int_{I_{a,\ b} }\ f (u )((a_2^2 -
(u - a_1)^2)^{n-2} \ (a_1 - u )\ d u ~ ,
\end{eqnarray*}
where $a=\exp(t_1H)$, $b=\exp(t_2H)$, $\delta(t)=(\sin
t)^{2(n-1)} \ \sin 2 t$,  $a_1 = \cos t_1 \cos t_2$, $a_2 = \sin
t_1 \sin t_2$ and
\begin{eqnarray*}
I_{a,b} = [\cos
 (t_1 + t_2) , \ \cos t_1 \cos t_2 - \sin t_1 \sin
t_2 \cos( \frac{\pi}{2(n-2)}) \ ].
\end{eqnarray*}

\end{thm}

\noindent {\bf Proof}-  Since both $f$ and $\lambda_a$ are
$K$-bi-invariant, $\lambda_a\star f(x)$ is $K$-bi-invariant and
therefore to compute $\mu_a\star f(x)$ we can assume that $x$ is
of the form $x= \exp(t'H)$.  Let $\{\theta_n\}$ be a sequence of
spherical functions converging weakly to the singular measure
$\lambda_a$ on the orbit $\mathcal{O}_a$.  Applying the polar
(Cartan) coordinate decomposition  for the convolution $\lambda_a
\star \check{f}(x)$ we have
\begin{eqnarray*}
\lambda_a \star \check{f}(x) &=& \int_{\mathcal{O}_a} f(y x^{-1})
dy
\\ &=& \lim_{n \rightarrow \infty} \int_G \theta_n(g) f(g x^{-1} ) dg
\\&=& c\lim_{n \rightarrow \infty} \int_K \int_K \int_A \theta_n(k_1 a' k_2 ) f ( k_1 a' k_2
x^{-1}) \delta(a') da' dk_1 dk_2 \\& =& c \lim_{n \rightarrow
\infty
} \int_K \int_A \theta_n(a') f(a' k x^{-1}) \delta(a') da' dk \\
&=& c \delta(t_1) \int_K f( a k x^{-1}) dk
\end{eqnarray*}
where
\begin{eqnarray*}
\delta(t) = \delta(\exp(tH))= (\sin t)^{2(n-1)} \sin 2t.
\end{eqnarray*}
Recall that $M$ is the centralizer group of $A$ in $K$ and it is
easy to verify that the function $g$ defined by $$ g(k) = \
f(\exp(t_1 H)k \exp(-t' H))$$ is an $M$- spherical function.
Applying the polar coordinates decomposition to the pair $(K, M)$
the above integral over $K$ becomes

\begin{eqnarray*}
\lambda_a \star \check{f}(x) &=& c \delta(t_1) \int_K f(\exp(t_1
H)k \exp(-t' H))d k \\&=& c\delta(t_1) \int_M \int_M \ \int_{A_0}
\ g(m_1 \ a \ m_2 )\delta_0 (a) d m_1 d a d m_2 \\ & = &  c
(\textrm{vol}.(M))^2 \delta(t_1) \int_{A_0} \ g(a) \delta_0(a) d
a ~,
\end{eqnarray*}
where $a = \exp( t H_0)$, $H_0=E_{23} - E_{32}$ is as above, $A_0$
is the corresponding real Cartan subgroup and $\delta_0$ is the
Jacobian of the polar coordinates corresponding to the pair
$(K,M)$. Thus we have

\begin{eqnarray*}
\lambda_a \star\check{f}(x)= c (\textrm{vol}.(M))^2 \delta(t_1)
\int_{Q_0} \ g(a) \delta_0( \exp( t H_0) d t
\end{eqnarray*}
where the polyhedron $Q_0$ is the interval $[0,
\frac{\pi}{2(n-2)}]$ in this case.  So the above convolution
integral becomes
\begin{eqnarray*}
\lambda_a \star \check{f}(x)= c \delta(t_1) (\textrm{vol}.(M))^2
\int_0^{\frac{\pi}{2(n-2)}}\ f ( \exp (t_1 H) \ \exp(t H_0) \exp
(-t' H)) \delta_0 ( \exp ( t H_0)) d t
\end{eqnarray*}
where $\delta(t_1)=(\sin t_1)^{2(n-1)} \ \sin 2 t_1$ and
$\delta_0(t) =(\sin t)^{2(n-2)} \sin 2t.$  We know that the
function $f$ is a spherical function and so it depends only to the
norm of the first entry of the product matrix
 $$\exp (t_1 H) \ \exp(t H_0) \exp (-t' H).
 $$
 After a  simple calculation we obtain $$ a_{11} = \cos t_1  \cos
 t'
  -   \cos t  \sin t_1  \sin t' .$$
 Set $a_1 =  \cos t_1  \cos t'$
 and $a_2 =  \sin t_1  \sin t'$ to obtain

\begin{eqnarray*}
 \lambda_a \star \check{f}(x) = c \delta(t_1) (\textrm{vol}.(M))^2 \int_0^{\frac{\pi}{2(n-2)}} \ f ( a_1  - a_2 \cos t
)\delta_0 ( \exp ( t H_0)) d t
\end{eqnarray*}
Next we compute the convolution $\lambda_a
\star \lambda_b(f)$ with $a=\exp(t_1 H)$ and $b=\exp(t_2 H)$.
Assume that $h(x)= \lambda_b \star \check{f}(x)$ then

\begin{eqnarray*}
(\lambda_a\star\lambda_b)(f)&=& (\lambda_a \star (\lambda_b
\star\check{f}))(e)\\&=&(\lambda_a \star h)(e) =\int_{O_a}
h(x)d\lambda_a (x)= h(a) \textrm{vol}(O_a),
\end{eqnarray*}
Applying polar coordinates on $G$ for the volume of the spherical
class $O_a$, $a=\exp(t_1 H)$, we obtain

\begin{eqnarray*}
 \textrm{vol} (O_a) = c \int_K \int_K  \delta(\exp(t_1 H))
 \ d k \ d k^{\prime} = c (\textrm{vol}(K))^2 (\sin ( t_1) )
^{2(n-1)} \ \sin 2 t_1.
\end{eqnarray*}
For $h(a)$ we have

\begin{eqnarray*}
 h(a)= \lambda_b \star \check{f}(a) = c \
\delta(t_2)(\textrm{vol}(M))^2 \ \int_0^{\frac{\pi}{2(n-2)}} \ f (
a_1  -  a_2 \cos t ) (\sin t)^{2(n-2)} \sin 2t \ d t.
\end{eqnarray*}
We make change of variable $ u = a_1 - a_2 \cos t$.  Then $u$ is
an increasing function on the interval $[0, \frac{ \pi}{2(n-2)}]$
and it maps this interval onto the interval
 $$I_{a,b} = [\cos
 (t_1 + t_2) , \cos t_1 \cos t_2 - \sin t_1 \sin
t_2 \cos( \frac{\pi}{2(n-2)}) ]
$$
Substituting the new variable $u$ and
simplifying the above integral we obtain

\begin{eqnarray*}
 h(a) =  \lambda_b \star \check{f}(a) = \frac{c \delta(t_2)(\textrm{vol}(M))^2}{a_2^{2(n-1)}} \
\int_{I_{a, b} }f (u ) (a_2^2 - (u - a_1)^2)^{n-2} (a_1 - u ) d u
\end{eqnarray*}
Finally for  $\lambda_a\star\lambda_b(f)$ we obtain

\begin{eqnarray*}
\lambda_a\star\lambda_b(f)&=& h(a) \textrm{vol}(O_a) \\& =&
\frac{(c\textrm{vol}(M)\textrm{vol}(K))^2}{a_2^{2(n-1)}}  \delta
(t_1) \delta (t_2) \int_{I_{a,\ b} } f (u )  (a_2^2 - (u -
 a_1)^2)^{n-2}(a_1 - u )d u
\end{eqnarray*}

which completes the proof of the theorem. $\blacksquare$\\

\begin{cor}
\label{cor:shaffaf13}Choosing matrices A and B according to the
(singular) invariant measures on the spherical classes ${\cal
O}_a$ and ${\cal O}_b$ respectively and normalized to be
probability measures, then the support of the distribution of the
product $A B$ is the interval $[\cos(t_2 + t_1), \cos t_1 \cos t_2
- \sin t_1 \sin t_2 \cos( \frac{\pi}{2(n-2)})]$ and its density
function is

\begin{eqnarray*}
\frac{2n-2}{(\sin t_1 \sin t_2 \sin
\frac{\pi}{2(n-2)})^{2n-2}}(a_2^2 - (u -  a_1)^2)^{n-2}  (a_1 - u
)~,
\end{eqnarray*}
 where $a_1 =
\cos t_1  \cos t_2$
and $a_2 =  \sin t_1  \sin t_2$.\\
\end{cor}

\begin{cor}
\label{cor:shaffaf51} With the notation and hypotheses of
Corollary \ref{cor:shaffaf13} the density function for the
convolution of probability measures $\lambda_a$ and $\lambda_b$
converges weakly to the singular invariant measure on the
spherical class through $\exp(t_1+t_2)H$ as $n\to\infty$.
\end{cor}

\noindent {\bf Proof} - Since $\lambda_a$ and $\lambda_b$ are
probability measures, so is $\lambda_a\star\lambda_b$.  The
support of this measure is the interval $I_{a,b}$, in the
appropriate coordinate system, which tends to the single point
$\cos (t_1+t_2)$ from which the required result follows.
$\blacksquare$

\section{ The Symmetric Pair $(SU(1,n), S(U(1) \times
U(n)))$.}

Let $G = Su(1,n)$ and $K = S(U(1) \times U(n))$, and
$\mathfrak{g}= \mathfrak{k}\oplus \mathfrak{p}$ be the
corresponding Cartan decomposition where $\mathfrak{g}$ and
$\mathfrak{k}$ are the Lie algebras of $G$ and $K$. Let  $
\mathfrak{a} \subset \mathfrak{p}$ be
\begin{eqnarray*}
\mathfrak{a} = \{ tH ~ ; ~~ t \in \mathbb{R} \}
\end{eqnarray*}
where $H= E_{12} + E_{21}$ and the restricted roots are given by
$\alpha( tH) = t$ and $2\alpha$ with multiplicities $2(n-1)$ and
$1$ respectively.  The centralizer of $\mathfrak{a}$ in
$\mathfrak{k}$ is $\mathfrak{m} = \mathfrak{s}(u(1) \times u(1)
\times u(n-1))$ and the corresponding Lie subgroup is $M = S(U(1)
\times U(1) \times U(n-1))$.  Therefore the pair $(K, M)$ in the
case of the non-compact symmetric pair $(SU(1,n), S(U(1) \times
U(n)))$ is same as in the case of the compact symmetric pair
$(SU(n+1), S(u(1) \times U(n)))$ treated in the preceding section.

\begin{thm}
\label{thm:shaffaf5} Let $\lambda_a$ and $\lambda_b$ be two
(singular) spherical measures concentrated on the $K$-spherical
classes ${\cal O}_a$ and ${\cal O}_b$ in the group $G= SU(1,n)$
respectively. Then $\lambda_a\star\lambda_b$ is absolutely
continuous relative to the Haar measure on $SU(1,n)$. It is a
spherical measure and for a continuous spherical function $f$ on
$SU(1,n)$ we have

\begin{eqnarray*}
\lambda_a\star\lambda_b(f) = (c\mathrm{vol}(M)\mathrm{vol}(K))^2
\delta (t_1) \delta (t_2) \int_{I_{a, b} } f (u )((a_2^2- (u -
a_1)^2)^{n-2}(a_1 - u ) d u ~ ,
\end{eqnarray*}
where $a=\exp(t_1H)$, $b=\exp(t_2H)$, $\delta(t)=(\sinh
t)^{2(n-1)} \sinh 2 t$, and
\begin{eqnarray*}
I_{a,b} = [\cosh
 (t_1 - t_2) ,  \cosh t_1 \cosh t_2 - \sinh t_1 \sinh
t_2 \cos( \frac{\pi}{2(n-2)}) \ ].
\end{eqnarray*}

\end{thm}

\noindent {\bf Proof}-  Since both $f$ and $\lambda_a$ are
$K$-bi-invariant, $\lambda_a\star f(x)$ is $K$-bi-invariant and
therefore to compute $\lambda_a\star f(x)$ we can assume that $x$
is of the form $x= \exp(t'H)$.  Let $\{\theta_n\}$ be a sequence
of spherical functions converging weakly to the singular measure
$\lambda_a$ on the orbit $\mathcal{O}_a$.  Applying the polar
coordinates (Cartan) decomposition  for the convolution $\lambda_a
\star \check{f}(x)$ we obtain
\begin{eqnarray*}
\lambda_a \star \check{f}(x) &=& \int_{\mathcal{O}_a} f(y x^{-1})
dy
\\ &=& \lim_{n \rightarrow \infty} \int_G \theta_n(g) f(g x^{-1} ) dg
\\&=& c\lim_{n \rightarrow \infty} \int_K \int_K \int_A \theta_n(k_1 a' k_2 ) f ( k_1 a' k_2
x^{-1}) \delta(a') da' dk_1 dk_2 \\& =& c \lim_{n \rightarrow
\infty
} \int_K \int_A \theta_n(a') f(a' k x^{-1}) \delta(a') da' dk \\
&=& c \delta(t_1) \int_K f( a k x^{-1}) dk
\end{eqnarray*}
where
\begin{eqnarray*}
\delta(t) = \delta(\exp(tH))= (\sinh t)^{2(n-1)} \sinh 2t.
\end{eqnarray*}
Recall that $M$ is the centralizer group of $A$ in $K$ and it is
easy to verify that the function $g$ defined by $$ g(k) =
f(\exp(t_1 H)k \exp(-t' H))$$ is an $M$-spherical function. Using
polar coordinates, as in the previous section, the above integral
over $K$ reduces to

\begin{eqnarray*}
\lambda_a \star \check{f}(x) &=& c \delta(t_1) \int_K f(\exp(t_1
H)k \exp(-t' H))d k \\&=& c\delta(t_1) \int_M \int_M \ \int_{A_0}
\ g(m_1 \ a \ m_2 )\delta_0 (a) d m_1 d a d m_2 \\ & = &  c
(\textrm{vol}.(M))^2 \delta(t_1) \int_{A_0} \ g(a) \delta_0(a) d
a ~,
\end{eqnarray*}
where $H_0=E_{12} - E_{21}$, $A_0=\{\exp (tH_0) \}$ is the real
Cartan subgroup and $\delta_0$ is as in the theorem
(\ref{thm:shaffaf1}).  Thus we have:

\begin{eqnarray*}
\lambda_a \star\check{f}(x)= c (\textrm{vol}.(M))^2 \delta(t_1)
\int_{Q_0} \ g(a) \delta_0( \exp( t H_0) d t ~,
\end{eqnarray*}
where the polyhedra $Q_0$ is the interval $[0,
\frac{\pi}{2(n-2)}]$.  Simplifying we obtain
\begin{eqnarray*}
\lambda_a \star \check{f}(x)= c \delta(t_1) (\textrm{vol}.(M))^2
\int_0^{\frac{\pi}{2(n-2)}}\ f ( \exp (t_1 H)  \exp(t H_0) \exp
(-t' H)) \delta_0 ( \exp ( t H_0)) d t
\end{eqnarray*}
where $\delta(t_1)=(\sinh t_1)^{2(n-1)} \ \sinh 2 t_1$ and
$\delta_0(t) =(\sin t)^{2(n-2)} \sin 2t.$  The function $f$ is
spherical and so it depends only on the norm of the first entry of
the product matrix
 $$\exp (t_1 H) \ \exp(t H_0) \exp (-t' H).
 $$
Now $$ a_{11} = \cosh t_1 \cosh
 t'
  - \cos t \sinh t_1 \sinh t'. $$
Set $a_1 =  \cosh t_1 \cosh t'$
 and $a_2 =  \sinh t_1 \sinh t'$ to obtain
\begin{eqnarray*}
 \lambda_a \star \check{f}(x) = c \delta(t_1) (\textrm{vol}.(M))^2 \int_0^{\frac{\pi}{2(n-2)}} \ f ( a_1 - a_2 \cos t
)\delta_0 ( \exp ( t H_0)) d t.
\end{eqnarray*}
Let $h(x)= \lambda_b \star \check{f}(x)$, then

\begin{eqnarray*}
(\lambda_a\star\lambda_b)(f)&=& (\lambda_a \star (\lambda_b
\star\check{f}))(e)\\&=&(\lambda_a \star h)(e) =\int_{O_a}
h(x)d\lambda_a (x)= h(a) \textrm{vol}(O_a).
\end{eqnarray*}
Using the decomposition $G=KAK$ we obtain
\begin{eqnarray*}
 \textrm{vol} (O_a) = c \int_K \int_K  \delta(\exp(t_1 H))
 \ d k \ d k^{\prime} = c (\textrm{vol}(K))^2 (\sinh ( t_1) )
^{2(n-1)} \ \sinh 2 t_1.
\end{eqnarray*}
Therefore
\begin{eqnarray*}
 h(a)= \lambda_b \star \check{f}(a) = c \
\delta(t_2)(\textrm{vol}(M))^2 \ \int_0^{\frac{\pi}{2(n-2)}} f (
a_1 -  a_2 \cos t ) (\sin t)^{2(n-2)} \sin 2t d t
\end{eqnarray*}
The change of variable $ u = a_1 - a_2 \cos t$ maps the interval
$[0, \frac{ \pi}{2(n-2)}]$ onto the interval
 $$I_{a, b} = [\cosh
 (t_1 - t_2) , \cosh t_1 \cosh t_2 - \sinh t_1  \sinh
t_2 \cos( \frac{\pi}{2(n-2)}) \ ]
$$
Therefore
\begin{eqnarray*}
 h(a) =  \lambda_b \star \check{f}(a) = \frac{ c \delta(t_2)(\textrm{vol}(M))^2}{a_2^{2(n-1)}}
\int_{I_{a, b} }\ f (u )(a_2^2 - (u  - a_1)^2)^{n-2} \
(a_1 - u ) d u
\end{eqnarray*}
Finally for  $\lambda_a\star\lambda_b(f)$ we obtain

\begin{eqnarray*}
\lambda_a\star\lambda_b(f)&=& h(a) \textrm{vol}(O_a) \\& =&
\frac{(c\textrm{vol}(M)\textrm{vol}(K))^2}{a_2^{2(n-1)}}  \delta
(t_1) \delta (t_2) \int_{I_{a,b} }\ f (u )(a_2^2 - (u -
 a_1)^2)^{n-2}(a_1 - u ) d u
\end{eqnarray*}
which completes the proof of the theorem. $\blacksquare$

\begin{cor}
\label{cor:shaffaf12}Choosing matrices A and B according to the
(singular) invariant measures on the spherical classes ${\cal
O}_a$ and ${\cal O}_b$ respectively and normalized to be
probability measures, then the support of the distribution of the
product $A B$ is the interval $[\cosh (t_2 - t_1), \cosh t_1 \
\cosh t_2 - \sinh t_1 \ \sinh t_2 \ \cos( \frac{\pi}{2(n-2)})  ]$
and its density function is
\begin{eqnarray*}
 \frac{2n-2}{(\sinh t_1 \sinh t_2 \sin \frac{\pi}{2(n-2)})^{2n-2}}(a_2^2 \ - (u \ -  \ a_1)^2)^{n-2} \ (a_1 - u )~,
\end{eqnarray*}
 where $a_1 =
\cosh t_1 \cosh t_2$ and $a_2 =  \sinh t_1 \sinh t_2$.
Furthermore
\begin{eqnarray*}
\lim_{n\to\infty} \lambda_a\star\lambda_b=\lambda_c,~~~~{\rm
weakly},
\end{eqnarray*}
where $\lambda_c$ is the singular invariant probability measure on
the spherical class through $\exp((t_1-t_2)H)$.
\end{cor}

\begin{rem}
\label{rem:shaffaf31}{\rm Note that $\exp(\pm tH)$ are in the same
spherical class and therefore $\exp((t_1-t_2)H)$ and
$\exp((t_2-t_1)H)$ are in the same spherical class.}
\end{rem}

\section{The symmetric pairs $(SO(n+1), SO(n))$, and $(SO(1,n)^\circ,SO(n))$, $n\ge 3$ }

The Lie algebra of the orthogonal group $G=SO(n+1)$ is the algebra
of skew symmetric matrices i.e.
$$\mathfrak{g}= \mathfrak{so}(n+1)= \{ A \ \in  M_{n+1} ( \mathbb{R} ) \ | \ A^t = - A \}
$$
For the Cartan decomposition of $\mathfrak{g}$ we have $$
\mathfrak{so}(n+1) = \mathfrak{so}(n) \oplus \mathfrak{p},$$ where
$\mathfrak{p}$ is the subspace spanned by the matrices of the form
\begin{eqnarray*}
\left[\begin{array}{ccccc}
0& & \xi_1 & \ldots & \xi_n  \\
\begin{array}{c}
\\ -\xi_1 \\
\vdots \\
\\ -\xi_n \\
\end{array} &  & & 0_{n \times n} & \\

\end{array}\right]
\end{eqnarray*}
A maximal abelian subspace of $\mathfrak{p}$ is:
$$ \mathfrak{a} = \{t H : t \in R \} $$ where $H$ is the matrix
\begin{eqnarray*} H = \left[\begin{array}{ccccc}
0& & \ 1 &  \ldots &  0 \\
\begin{array}{c}
\ -1 \\
\ 0 \\
\vdots \\
\ 0 \\
\end{array} &  & & 0_{n \times n} & \\

\end{array}\right]
\end{eqnarray*}

\noindent The centralizer $\mathfrak{m}$ of $\mathfrak{a}$ in
$\mathfrak{k}= \mathfrak{so}(n)$ is exactly the Lie algebra
$\mathfrak{so}(n-1)$. It is straightforward that the centralizer
subgroup $M$ is connected and therefore we have $M=SO(n-1)$ and
 $$ K  / M =
SO(n) / SO(n-1) \cong S^{n-1}. $$ By straightforward calculation
the eigenvalues of the operator
$$\textrm{ad} H : \mathfrak{so}(n) \longrightarrow\mathfrak{so}(n)$$
are $\pm i$ with multiplicity $n-1$. Thus for $(SO(n+1), SO(n))$
we have one positive restricted root $\alpha(tH) = it$ whose
multiplicity is $m_\alpha = n-1$.  Hence
$$ \delta(\exp(t H)) = \prod_{ \alpha \in \  \Sigma^+}  | \sin
\alpha (i t H)|^ {m_ {\alpha}} = (\sin ( t) ) ^{n-1}
$$
For the pair $(K,M) = (SO(n), SO(n-1))$ the  maximal abelian
subspace $\mathfrak{a}_0$ is
$$ \mathfrak{a}_0 = \{ t H_0 \ | \ t \in \mathbb{R} \} \ \ , $$
where $H_0$ is
 $$H_0= \left[\begin{array}{c|c|ccc} 0&0&0&\cdots&0
\\\hline 0&0& 1 & \cdots & 0 \\\hline
0& -1   &&&\\\vdots&\vdots&&0_{n-3}&\\ 0 & 0 &&&\\
\end{array}\right]
$$
Therefore the corresponding Jacobian for the pair $(K,M)$ is
$$
\delta_0(\exp(t H_0)) = \prod_{ \alpha \in \ \Sigma^+} | \sin
\alpha (i t H_0)|^ {m_ {\alpha}} = (\sin ( t) ) ^{n-2}.
$$
Since $SO(n+1)$ is not simply connected we work with the double
cover ${\rm Spin}(n+1)$ which is simply connected, and let
$$\pi :{\rm Spin(n+1)} \longrightarrow SO(n+1).$$
Since $ \mathfrak{spin}(n+1) = \mathfrak{so}(n+1)$ we have $
\mathfrak{spin}(n+1) = \mathfrak{spin}(n) \oplus \mathfrak{p}$.
Note that we identify ${\rm Spin}(n-1)$ with the pre-image of the
subgroup $K = SO(n-1)$ under the covering map $\pi$.  We denote
this subgroup by
$$\widetilde{K} = {\rm Spin}(n) = \pi^{-1} (SO(n) ).$$
The computation of the restricted positive root and its
multiplicity is the same as in the case of $\mathfrak{so}(n+1)$.
Now Theorem (\ref{thm:shaffaf1}) is applicable to the symmetric
pair $({\rm Spin}(n+1), {\rm Spin}(n) )$, but note that a function
on the group $SO(n+1)$ can be considered as a function on ${\rm
Spin}(n)$ and its integral over the group $G = {\rm Spin}(n+1)$ is
equal to twice its integral over $G= SO(n+1)$.

\begin{thm}
\label{thm:shaffaf23} Let $\lambda_a$ and $\lambda_b$ be two
(singular) spherical measures concentrated on the $K$-spherical
classes ${\cal O}_a$ and ${\cal O}_b$ in the group $G= SO(n)$
respectively. Then $\lambda_a\star\lambda_b$ is absolutely
continuous relative to the Haar measure on $SO(n)$. It is a
spherical measure and for a continuous spherical function $f$ on
$SO(n)$ we have

\begin{eqnarray*}
\lambda_a\star\lambda_b(f) = (c\mathrm{vol}(M)\mathrm{vol}(K))^2
\delta (t_1) \delta (t_2) \int_{I_{a,b} }\ f (u )((a_2^2 - (u -
a_1)^2)^{\frac{n-4}{2}} d u ~ ,
\end{eqnarray*}
where $a=\exp(t_1H)$, $b=\exp(t_2H)$, $\delta(t)=(\sin t)^{n-2}$,
and
\begin{eqnarray*}
I_{a,b} = [\cos
 (t_1 + t_2) , \cos t_1  \cos t_2 - \sin t_1  \sin
t_2 \cos( \frac{\pi}{n-3)})  ].
\end{eqnarray*}

\end{thm}

\noindent {\bf Proof}-  Since both $f$ and $\lambda_a$ are
$K$-bi-invariant, $\lambda_a\star f(x)$ is $K$-bi-invariant and
therefore to compute $\lambda_a\star f(x)$ we can assume that $x$
is of the form $x= \exp(t'H)$.  Let $\{\theta_n\}$ be a sequence
of spherical functions converging weakly to the singular measure
$\lambda_a$ on the orbit $\mathcal{O}_a$.  Applying the polar
coordinates (Cartan) decomposition  for the convolution $\lambda_a
\star \check{f}(x)$ we have
\begin{eqnarray*}
\lambda_a \star \check{f}(x)
&=& \lim_{n \rightarrow \infty} \int_G \theta_n(g) f(g x^{-1} ) dg
\\&=& \frac{1}{2} \ \textrm{lim}
\int_{\widetilde{G}} \theta_n(g) f(g x^{-1}) d g \\& =&\frac{
c}{2}\lim_{n \rightarrow \infty} \int_{ \widetilde{K}} \int_{
\widetilde{K}} \int_{\widetilde{A}} \theta_n(k_1 a' k_2 ) f ( k_1
a' k_2 x^{-1}) \delta(a') da' dk_1 dk_2 \\& =& \frac{c}{2} \lim_{n
\rightarrow \infty
} \int_{\widetilde{K}} \int_{\widetilde{A}} \theta_n(a') f(a' k x^{-1}) \delta(a') da' dk \\
&=&\frac{c}{2} \delta(t_1) \int_{\widetilde{K}} f( a k x^{-1}) dk
\\&=& c \delta(t_1) \int_K f( a k x^{-1}) dk
\end{eqnarray*}
where
\begin{eqnarray*}
\delta(t) = \delta(\exp(tH))= (\sin t)^{n-2}.
\end{eqnarray*}
Recall that $M$ is the centralizer of $A$ in $K$ and that the
function $g$ defined by $$ g(k) = \ f(\exp(t_1 H)k \exp(-t' H))$$
is an $M$-spherical function. Applying the polar coordinates
decomposition to the pair $(K, M)$ the above integral becomes
\begin{eqnarray*}
\lambda_a \star \check{f}(x) &=& c \delta(t_1) \int_K f(\exp(t_1
H)k \exp(-t' H))d k \\&=& c\delta(t_1) \int_M \int_M \ \int_{A_0}
\ g(m_1 \ a \ m_2 )\delta_0 (a) d m_1 d a d m_2 \\ & = &  c
(\textrm{vol}.(M))^2 \delta(t_1) \int_{A_0} \ g(a) \delta_0(a) d
a ~,
\end{eqnarray*}
where $a = \exp( t H_0)$, where $H_0=E_{23} - E_{32}$ is as above,
$A_0$ is the real Cartan subgroup and $\delta_0$ is the Jacobian
of the polar coordinates corresponding to the pair $(K,M)$. Thus
we have:

\begin{eqnarray*}
\lambda_a \star\check{f}(x)= c (\textrm{vol}.(M))^2 \delta(t_1)
\int_{Q_0} \ g(a) \delta_0( \exp( t H_0) d t ~,
\end{eqnarray*}
where the polyhedra $Q_0$ is the interval $[0, \frac{\pi}{n-2}]$.
Therefore
\begin{eqnarray*}
\lambda_a \star \check{f}(x)= c \delta(t_1) (\textrm{vol}.(M))^2
\int_0^{\frac{\pi}{n-2}}\ f ( \exp (t_1 H) \ \exp(t H_0) \exp (-t'
H)) \delta_0 ( \exp ( t H_0)) d t
\end{eqnarray*}
where $\delta(t_1)=(\sin t_1)^{n-1}$ and $\delta_0(t) =(\sin
t)^{n-2}.$  The function $f$ is spherical and so it depends only
on the norm of the first entry of the product matrix $\exp (t_1 H)
 \exp(t H_0) \exp (-t' H),$ and is given by $ a_{11} = \cos t_1
\cos t'  -  \cos t  \sin t_1 \ \sin t'.$ Set  $a_1 = \cos t_1 \cos
t'$ and $a_2 =  \sin t_1 \sin t'$ to obtain

\begin{eqnarray*}
 \lambda_a \star \check{f}(x) = c \delta(t_1) (\textrm{vol}.(M))^2 \int_0^{\frac{\pi}{n-3}} f ( a_1  -  a_2 \cos t
)\delta_0 ( \exp ( t H_0)) d t
\end{eqnarray*}
Let $h(x)= \lambda_b \star \check{f}(x)$ with $a=\exp(t_1 H)$ and
$b=\exp(t_2 H)$, then

\begin{eqnarray*}
(\lambda_a\star\lambda_b)(f)&=& (\lambda_a \star (\lambda_b
\star\check{f}))(e)\\&=&(\lambda_a \star h)(e)\\ & =&\int_{O_a}
h(x)d\lambda_a (x)= h(a) \textrm{vol}(O_a),
\end{eqnarray*}
Using the Cartan decomposition we obtain

\begin{eqnarray*}
 \textrm{vol} (O_a) =\frac {c}{2} \int_{\widetilde{K}} \int_{\widetilde{K}}  \delta(\exp(t_1 H))
 \ d k \ d k^{\prime} = 2 c (\textrm{vol}(K))^2 (\sin ( t_1) )
^{n-1}
\end{eqnarray*}
Now
\begin{eqnarray*}
 h(a)= \lambda_b \star \check{f}(a) = c \
\delta(t_2)(\textrm{vol}(M))^2 \ \int_0^{\frac{\pi}{n-2}} \ f (
a_1  -  a_2 \cos t ) (\sin t)^{n-2} \ d t
\end{eqnarray*}
The change of variable $ u = a_1  -  a_2 \cos t$ maps the interval
$[0, \frac{ \pi}{n-2}]$ onto the interval
 $$I_{a,b} = [\cos
 (t_1 + t_2) , \cos t_1 \cos t_2 - \sin t_1 \sin
t_2 \cos( \frac{\pi}{n-2}) \ ]
$$
The expression for $h(a)$ becomes

\begin{eqnarray*}
 h(a) =  \lambda_b \star \check{f}(a) = \frac{c \delta(t_2)(\textrm{vol}(M))^2}{a_2^{n-2}} \
\int_{I_{a,b} }  f (u ) (a_2^2  - (u - a_1)^2)^{\frac{n-3}{2}} d u
\end{eqnarray*}

Finally for  $\lambda_a\star\lambda_b(f)$ we obtain

\begin{eqnarray*}
\lambda_a\star\lambda_b(f)&=& h(a) \textrm{vol}(O_a) \\& =&
\frac{(c\textrm{vol}(M)\textrm{vol}(K))^2  \delta (t_1) \delta
(t_2)}{a_2^{n-2}} \int_{I_{a,b} } \  f (u ) (a_2^2  - (u -
a_1)^2)^{\frac{n-3}{2}} d u
\end{eqnarray*}
which completes the proof of the theorem. $\blacksquare$

\begin{cor}
\label{cor:shaffaf26}Choosing matrices A and B according to the
$($singular$)$ invariant measures on the spherical classes ${\cal
O}_a$ and ${\cal O}_b$ respectively and normalized to be probability
measure, then the support of the distribution of the product $A B$
is the interval
\[[\cos (t_2 + t_1), \cos t_1 \cos t_2 - \sin t_1 \sin t_2 \cos(
\frac{\pi}{n-2})]\] and its density function is
\begin{eqnarray*}
\frac{(a_2^2  - (u -  a_1)^2)^{\frac{n-3}{2}}}{a_2^{n-2}
\int_0^{\frac{\pi}{n-2}} (\sin t)^{n-2} dt} ~,
\end{eqnarray*}
where $a_1=\cos t_1\cos t_2$ and $a_2=\sin t_1\sin t_2$.
Furthermore
\begin{eqnarray*}
\lim_{n\to\infty} \lambda_a\star\lambda_b=\lambda_c,~~~~{\rm
weakly},
\end{eqnarray*}
where $\lambda_c$ is the singular invariant probability measure on
the spherical class through $\exp((t_1+t_2)H)$.
\end{cor}

For the symmetric pair $(SO(1,n)^\circ,SO(n))$ the calculations
are similar and therefore are not repeated.  We obtain

\begin{thm}
\label{thm:shaffaf24} Let $\lambda_a$ and $\lambda_b$ be two
(singular) spherical measures concentrated on the $K$-spherical
classes ${\cal O}_a$ and ${\cal O}_b$ in the group $G= SO(1,n)$
respectively. Then $\lambda_a\star\lambda_b$ is absolutely
continuous relative to the Haar measure on $SO(1,n)$. It is a
spherical measure and for a continuous spherical function $f$ on
$SO(1,n)$ we have

\begin{eqnarray*}
\lambda_a\star\lambda_b(f) = (c\rm{vol}(M)\rm{vol}(K))^2 \delta
(t_1) \delta (t_2) \int_{I_{a,b} }\ f (u )((a_2^2 - (u -
a_1)^2)^{\frac{n-3}{2}} d u ~ ,
\end{eqnarray*}
where $a=\exp(t_1H)$, $b=\exp(t_2H)$, $\delta(t)=(\sinh t)^{n-1}$,
and
\begin{eqnarray*}
I_{a,b} = [\cosh
 (t_1 - t_2) , \cosh t_1  \cosh t_2 - \sinh t_1  \sinh
t_2 \cos( \frac{\pi}{n-2})  ].
\end{eqnarray*}
\end{thm}

\begin{cor}
\label{cor:shaffaf27}Choosing matrices A and B according to the
$($singular$)$ invariant measures on the spherical classes ${\cal
O}_a$ and ${\cal O}_b$ respectively and normalized to be probability
measure, then the support of the distribution of the product $A B$
is the interval
\[[\cosh (t_1 - t_2), \cosh t_1 \cosh t_2 - \sinh t_1 \sinh t_2 \cos(
\frac{\pi}{n-2})]\] and its density function is
\begin{eqnarray*}
\frac{(a_2^2  - (u -  a_1)^2)^{\frac{n-3}{2}}}{a_2^{n-2}
\int_0^{\frac{\pi}{n-2}} (\sin t)^{n-2} dt} ~,
\end{eqnarray*}
where $a_1=\cosh t_1\cosh t_2$ and $a_2=\sinh t_1\sinh t_2$.
Furthermore
\begin{eqnarray*}
\lim_{n\to\infty} \lambda_a\star\lambda_b=\lambda_c,~~~~{\rm
weakly},
\end{eqnarray*}
where $\lambda_c$ is the singular invariant probability measure on
the spherical class through $\exp((t_1-t_2)H)$.
\end{cor}

\section{The symmetric pairs $(Sp(n+1), Sp(1)\times Sp(n))$, and $(Sp(1,n),Sp(1)\times
Sp(n))$}

Let $G= Sp(n+1)$ and $K = Sp(1) \times Sp(n)$, then $G$ is simply
connected, $K$ connected and $G/K$ is the quaternionic projective
space. The Lie algebra of $G$ is $\mathfrak{g}= sp(n+1)$, the
space of complex matrices $X$ satisfying $ JX + X^t J = 0$.  If we
write $X$ in the form
$$\left[\begin{array}{c|c}
X_1&X_2\\\hline
X_3&X_4
\end{array}\right] $$
where $X_1, X_2, X_3, X_4$ are  matrices of degree $n+1$, the
condition $ JX + X^t J = 0$ gives
\begin{eqnarray*}
X_4 = - X_1^t  ~~~ X_3 = X_3^t ~~~ X_2 = X_2^t
\end{eqnarray*}
and the Lie algebra of the subgroup $K$ is
 \begin{eqnarray*}
\mathfrak{k}=
\left\{\left[\begin{array}{cc|cc}
x_{11}&x_{12}&0&0\\
-\overline{x_{12}}&\overline{x_{11}}&0&0\\\hline
0&0&Y_{11}&Y_{12}\\
0&0&-\overline{Y_{12}}&\overline{Y_{11}}
\end{array}\right]; \
\textrm{$x_{ij} \in \mathbb{C}$ and $Y_{11} \in \mathfrak{u}(n)$, $Y_{12}$ is $ n \times n$ symmetric}\right\}.
\end{eqnarray*}
Let $\mathfrak{g} = \mathfrak{k} \oplus \mathfrak{p}$ be the Cartan
decomposition then $\mathfrak{p}$ is the subspace of matrices of the
form
\begin{eqnarray*}
\left[\begin{array}{cc|cc}
0&0&Z_{13}&-\overline{Z_{14}}\\
0&0&-Z_{14}&-\overline{Z_{13}}\\\hline
-\overline{Z_{13}}^t&\overline{Z_{14}}^t&0&0\\
Z_{14}^t&Z_{13}^t&0&0
\end{array}\right] \ ~,
\end{eqnarray*}
where $Z_{ij}$ are $1 \times n$  complex matrices.
A maximal abelian subspace $ \mathfrak{a} \subset \mathfrak{p}$ is $
\mathfrak{a} = \{t H ~;~~ t \in \mathbb{R} \}$ where $H$ is the matrix

 \begin{eqnarray*}
  H = E_{32} + E_{41}-E_{23}-E_{14}.
\end{eqnarray*}
The restricted roots are given by $\alpha( tH) = it$ and $2\alpha$
with multiplicities $8(n-1)$ and $2$ respectively.  The
centralizer of $\mathfrak{a}$ in $\mathfrak{k}$ is $\mathfrak{m} =
\mathfrak{sp}(1) \times \mathfrak{sp}(1) \times
\mathfrak{sp}(n-1)$, the corresponding Lie subgroup $M$ is
connected and $M = Sp(1) \times Sp(1) \times Sp(n-1)$.

\begin{thm}
\label{thm:shaffaf22} Let $\lambda_a$ and $\lambda_b$ be two
(singular) spherical measures concentrated on the $K$-spherical
classes ${\cal O}_a$ and ${\cal O}_b$ in the group $G= Sp(n+1)$
respectively. Then $\lambda_a\star\lambda_b$ is absolutely
continuous relative to the Haar measure on $Sp(n+1)$. It is a
spherical measure and for a continuous spherical function $f$ on
$Sp(n+1)$ we have

\begin{eqnarray*}
\lambda_a\star\lambda_b(f) =
\frac{(c\rm{vol}(M)\rm{vol}(K))^2}{a_2^{8n-12}}  \delta (t_1)
\delta (t_2) \int_{I_{a, b} } f (u )(a_2^2 - (u - a_1)^2)^{4n-
\frac{15}{2}} (a_1 - u )^2 d u
 ~ ,
\end{eqnarray*}
where $c$ is a constant, $a=\exp(t_1H)$, $b=\exp(t_2H)$,
$\delta(t)=(\sin t)^{8(n-1)} \ (\sin 2 t)^2$,  $a_1 = \cos t_1
\cos t_2$, $a_2 = \sin t_1 \sin t_2$ and
\begin{eqnarray*}
I_{a,b} = [\cos
 (t_1 + t_2) , \ \cos t_1 \cos t_2 - \sin t_1 \sin
t_2 \cos( \frac{\pi}{8(n-2)})  ].
\end{eqnarray*}

\end{thm}

\noindent {\bf Proof}-  Since both $f$ and $\lambda_a$ are
$K$-bi-invariant, $\lambda_a\star f(x)$ is $K$-bi-invariant and
therefore to compute $\mu_a\star f(x)$ we can assume that $x$ is of
the form $x= \exp(t'H)$.  Let $\{\theta_n\}$ be a sequence of
spherical functions converging weakly to the singular measure
$\lambda_a$ on the orbit $\mathcal{O}_a$.  Applying the polar
(Cartan) coordinate decomposition  for the convolution $\lambda_a
\star \check{f}(x)$ we have
\begin{eqnarray*}
\lambda_a \star \check{f}(x) &=& \int_{\mathcal{O}_a} f(y x^{-1}) dy
\\ &=& \lim_{n \rightarrow \infty} \int_G \theta_n(g) f(g x^{-1} ) dg
\\&=& c\lim_{n \rightarrow \infty} \int_K \int_K \int_A \theta_n(k_1 a' k_2 ) f ( k_1 a' k_2
x^{-1}) \delta(a') da' dk_1 dk_2 \\& =& c \lim_{n \rightarrow \infty
} \int_K \int_A \theta_n(a') f(a' k x^{-1}) \delta(a') da' dk \\
&=& c \delta(t_1) \int_K f( a k x^{-1}) dk
\end{eqnarray*}
where
\begin{eqnarray*}
\delta(t) = \delta(\exp(tH))= (\sin t)^{8(n-1)} (\sin 2t)^2.
\end{eqnarray*}
The function $g$ defined by $$ g(k) = \ f(\exp(t_1 H)k \exp(-t'
H))$$ is an $M$- spherical function. Applying the polar
coordinates decomposition to the pair $(K,M)$, the above integral
over $K$ becomes

\begin{eqnarray*}
\lambda_a \star \check{f}(x) &=& c \delta(t_1) \int_K f(\exp(t_1 H)k
\exp(-t' H))d k \\&=& c\delta(t_1) \int_M \int_M \ \int_{A_0} \
g(m_1 \ a \ m_2 )\delta_0 (a) d m_1 d a d m_2 \\ & = &  c
(\textrm{vol}.(M))^2 \delta(t_1) \int_{A_0} \ g(a) \delta_0(a) d a
~,
\end{eqnarray*}
where $a = \exp( t H_0)$, $H_0=E_{74} + E_{83}-E_{56} -E_{65}$ is as
above, $A_0$ is the corresponding real Cartan subgroup and
$\delta_0$ is the Jacobian of the polar coordinates corresponding to
the pair $(K,M)$. Thus we have

\begin{eqnarray*}
\lambda_a \star\check{f}(x)= c (\textrm{vol}.(M))^2 \delta(t_1)
\int_{Q_0} \ g(a) \delta_0( \exp( t H_0) d t
\end{eqnarray*}
where the polyhedron $Q_0$ is the interval $[0, \frac{\pi}{8(n-2)}]$
in this case.  So the above convolution integral becomes
\begin{eqnarray*}
\lambda_a \star \check{f}(x)= c \delta(t_1) (\textrm{vol}.(M))^2
\int_0^{\frac{\pi}{8(n-2)}}\ f ( \exp (t_1 H) \ \exp(t H_0) \exp
(-t' H)) \delta_0 ( \exp ( t H_0)) d t
\end{eqnarray*}
where $\delta(t_1)=(\sin t_1)^{8(n-1)} \ (\sin 2 t_1)^2$ and
$\delta_0(t) =(\sin t)^{8(n-2)} (\sin 2t)^2.$  We know that the
function $f$ is a spherical function and so it depends only on the
norm of the first entry $a_{11}$ of the product matrix
 $\exp (t_1 H) \ \exp(t H_0) \exp (-t' H)$, and
 after a  simple calculation we obtain $$ a_{11} = \cos t_1  \cos
 t'
  -   \cos t  \sin t_1  \sin t' .$$
 Set $a_1 =  \cos t_1  \cos t'$
 and $a_2 =  \sin t_1  \sin t'$ to obtain

\begin{eqnarray*}
 \lambda_a \star \check{f}(x) = c \delta(t_1) (\textrm{vol}.(M))^2 \int_0^{\frac{\pi}{8(n-2)}} \ f ( a_1  - a_2 \cos t
)\delta_0 ( \exp ( t H_0)) d t
\end{eqnarray*}
Next we compute the convolution $\lambda_a \star \lambda_b(f)$ with
$a=\exp(t_1 H)$ and $b=\exp(t_2 H)$. Assume that $h(x)= \lambda_b
\star \check{f}(x)$ then

\begin{eqnarray*}
(\lambda_a\star\lambda_b)(f)&=& (\lambda_a \star (\lambda_b
\star\check{f}))(e)\\&=&(\lambda_a \star h)(e) =\int_{O_a}
h(x)d\lambda_a (x)= h(a) \textrm{vol}(O_a),
\end{eqnarray*}
Applying polar coordinates on $G$ for the volume of the spherical
class $O_a$, $a=\exp(t_1 H)$, we obtain

\begin{eqnarray*}
 \textrm{vol} (O_a) = c \int_K \int_K  \delta(\exp(t_1 H))
 \ d k \ d k^{\prime} = c (\textrm{vol}(K))^2 (\sin ( t_1) )
^{8(n-1)} \ (\sin 2 t_1)^2.
\end{eqnarray*}
For $h(a)$ we have

\begin{eqnarray*}
 h(a)= \lambda_b \star \check{f}(a) = c \
\delta(t_2)(\textrm{vol}(M))^2 \ \int_0^{\frac{\pi}{8(n-2)}} \ f (
a_1  -  a_2 \cos t ) (\sin t)^{8(n-2)} (\sin 2t)^2 \ d t.
\end{eqnarray*}
We make change of variable $ u = a_1 - a_2 \cos t$.  Then $u$ is an
increasing function on the interval $[0, \frac{ \pi}{8(n-2)}]$ and
it maps this interval onto the interval
 $$I_{a,b} = [\cos
 (t_1 + t_2) , \cos t_1 \cos t_2 - \sin t_1 \sin
t_2 \cos( \frac{\pi}{8(n-2)}) ]
$$
Substituting the new variable $u$ and simplifying the above integral
we obtain

\begin{eqnarray*}
 h(a) =  \lambda_b \star \check{f}(a) = \frac{c \delta(t_2)(\textrm{vol}(M))^2}{a_2^{8n-12}} \
\int_{I_{a, b} }f (u ) (a_2^2 - (u - a_1)^2)^{4n- \frac{15}{2}} (a_1
- u )^2 d u
\end{eqnarray*}
Finally for  $\lambda_a\star\lambda_b(f)$ we obtain

\begin{eqnarray*}
\lambda_a\star\lambda_b(f)&=& h(a) \textrm{vol}(O_a) \\& =&
\frac{(c\textrm{vol}(M)\textrm{vol}(K))^2}{a_2^{8n-12}}  \delta
(t_1) \delta (t_2) \int_{I_{a, b} } f (u )(a_2^2 - (u -
a_1)^2)^{4n- \frac{15}{2}} (a_1 - u )^2 d u
\end{eqnarray*}

which completes the proof of the theorem. $\blacksquare$\\

\begin{cor}
\label{cor:shaffaf23}Choosing matrices A and B according to the
(singular) invariant measures on the spherical classes ${\cal O}_a$
and ${\cal O}_b$ respectively and normalized to be probability
measures, then the support of the distribution of the product $A B$
is the interval $[\cos(t_2 + t_1), \cos t_1 \cos t_2 - \sin t_1 \sin
t_2 \cos( \frac{\pi}{8(n-2)})]$ and its density function is

\begin{eqnarray*}
\frac{1}{a_2^{8n-12}\int_0^{\frac{\pi}{8(n-2)}}(\sin t)^{8n-14}(\cos
t)^2 dt}(a_2^2 - (u - a_1)^2)^{4n- \frac{15}{2}} (a_1 - u )^2~,
\end{eqnarray*}
 where $a_1 =
\cos t_1  \cos t_2$ and $a_2 =  \sin t_1  \sin t_2$. Furthermore
\begin{eqnarray*}
\lim_{n\to\infty} \lambda_a\star\lambda_b=\lambda_c,~~~~{\rm
weakly},
\end{eqnarray*}
where $\lambda_c$ is the singular invariant probability measure on
the spherical class through $\exp((t_1+t_2)H)$.
\end{cor}
For the symmetric pair $(Sp(1,n),Sp(1) \times Sp(n))$ the
calculations are similar and therefore are not repeated.  We obtain

\begin{thm}
\label{thm:shaffaf222} Let $\lambda_a$ and $\lambda_b$ be two
(singular) spherical measures concentrated on the $K$-spherical
classes ${\cal O}_a$ and ${\cal O}_b$ in the group $G= Sp(1,n)$
respectively. Then $\lambda_a\star\lambda_b$ is absolutely
continuous relative to the Haar measure on $Sp(1,n)$. It is a
spherical measure and for a continuous spherical function $f$ on
$Sp(1,n)$ we have

\begin{eqnarray*}
\lambda_a\star\lambda_b(f) =
\frac{(c\rm{vol}(M)\rm{vol}(K))^2}{a_2^{8n-12}}  \delta (t_1) \delta
(t_2) \int_{I_{a,\ b} } f (u )(a_2^2 - (u - a_1)^2)^{4n-
\frac{15}{2}} (a_1 - u )^2 d u
 ~ ,
\end{eqnarray*}
where $a=\exp(t_1H)$, $b=\exp(t_2H)$, $\delta(t)=(\sinh
t)^{8(n-1)}(\sinh 2t)^2$, and $a_1 = \cosh t_1 \cosh t_2$, $a_2 =
\sinh t_1 \sinh t_2$
\begin{eqnarray*}
I_{a,b} = [\cosh
 (t_1 - t_2) , \cosh t_1  \cosh t_2 - \sinh t_1  \sinh
t_2 \cos( \frac{\pi}{8(n-2)})  ].
\end{eqnarray*}
\end{thm}

\begin{cor} \label{cor:shaffaf223}Choosing
matrices A and B according to the (singular) invariant measures on
the spherical classes ${\cal O}_a$ and ${\cal O}_b$ respectively and
normalized to be probability measures, then the support of the
distribution of the product $A B$ is the interval $[\cosh (t_2
-t_1), \cosh t_1 \cosh t_2 - \sinh t_1 \sinh t_2 \cos(
\frac{\pi}{8(n-2)})]$ and its density function is

\begin{eqnarray*}
\frac{1}{a_2^{8n-12}\int_0^{\frac{\pi}{8(n-2)}}(\sin t)^{8n-14}(\cos
t)^2 dt}(a_2^2 - (u - a_1)^2)^{4n- \frac{15}{2}} (a_1 - u )^2~,
\end{eqnarray*}
 where $a_1 =
\cosh t_1  \cosh t_2$ and $a_2 =  \sinh t_1  \sinh t_2$.
Furthermore
\begin{eqnarray*}
\lim_{n\to\infty} \lambda_a\star\lambda_b=\lambda_c,~~~~{\rm
weakly},
\end{eqnarray*}
where $\lambda_c$ is the singular invariant probability measure on
the spherical class through $\exp((t_1-t_2)H)$.
\end{cor}

\section{Convergence to Haar Measure}
It was noted that $\lambda_a\star\lambda_b$ behaves approximately
as in the abelian case when $n\to\infty$.  In this section we
determine the rate of convergence of
$(\lambda_a\star\lambda_b)^{l(n)}$ to the Haar measure as
$n\to\infty$. More precisely we prove
\begin{thm}
\label{thm:jshaffaf2} Let $\lambda_a$ denote the invariant measure
on the spherical class ${\mathcal O}_a$.  Then
$(\lambda_a\star\lambda_b)^{l(n)}$ converges to the Haar measure
on $SU(n)$ as $n\to\infty$ if $l(n)\ge c\log n$ where $c$ is a
constant depending on the choice of the spherical classes
${\mathcal O}_a$ and ${\mathcal O}_b$.
\end{thm}

To interpret this theorem let ${\mathcal O}_a$ and ${\mathcal
O}_b$ denote spherical classes where $a=\exp(t_1H)$ and
$b=\exp(t_2H)$. We had observed that as $n\to\infty$ the Product
${\mathcal O}_a.{\mathcal O}_b$ converges to the spherical measure
concentrated on the spherical class passing through
$\exp(t_1+t_2)H$.  The measure $(\lambda_a\star\lambda_b)^{l(n)}$
represents the empirical measure of products

$$ A_1 B_1 A_2 B_2 \ldots A_{l(n)}B_{l(n)}$$

\noindent where $A_i$'s are chosen randomly on ${\mathcal O}_a$
and similarly for $B_j$'s.  Theorem \ref{thm:jshaffaf2} asserts
that for $l(n)$ of the stated form this empirical measure
converges weakly to the Haar measure on $G=SU(n)$ as $n\to\infty$.

The proof of this theorem requires Schur-Weyl theory on the
representations of $SU(n)$ (see [B] or [W] for an account of
Schur-Weyl theory). To fix notation we recall the relevant facts.
Let $T$ denote a Young diagram, i.e., a graphical representation
of a partition $m=m_1+\ldots +m_k$, $k\le n$, with $m_i\ge
m_{i+1}$.  A Young diagram $T$ filled with integers $1,2,\ldots,m$
is denoted by $\{T\}$ and called a Young tableau. A standard Young
tableau is one such that the integers are strictly increasing
along rows and columns.   We fix the enumeration of the boxes in a
Young diagram by starting at the upper left corner and moving
along columns consecutively. With this enumeration of the squares
in a Young diagram $T$, two subgroups of the symmetric group
${\cal S}_m$ specified, namely, the subgroup $H=H_T$ consisting of
all permutations preserving the rows, and $H^\prime=H^\prime_T$
consisting of all permutations preserving the columns.  Let ${\bf
Z}[{\cal S}_m]$ be the {\it integral group algebra} of the
symmetric group, i.e., formal linear combinations with integers
coefficients and multiplication inherited from the group law in
${\cal S}_m$. Define the {\it Young symmetrizer}
${\sf{C}}={\sf{C}}_T\in {\bf Z}[{\cal S}_m]$ as
\begin{eqnarray*}
{\sf{C}}={\sf{C}}_T=(\sum_{\tau\in H^\prime}\epsilon_\tau \tau)
(\sum_{\sigma\in H}\sigma)= \sum_{\sigma\in H,\tau\in H^\prime}
\epsilon_\tau \tau\sigma.
\end{eqnarray*}

The symmetric group ${\cal S}_m$ and therefore its group algebra
${\bf Z}[{\cal S}_m]$ act on the tensor space $T^m(V)$.  In fact,
given a tensor $v_{i_1}\otimes\cdots\otimes v_{i_m}$, $v_{i_j}\in
V$, and $\sigma\in {\cal S}_m$, the action of $\sigma$ is given by
\begin{eqnarray*}
v_{i_1}\otimes\cdots\otimes v_{i_m}
\buildrel\sigma\over\longrightarrow
v_{i_{\sigma(1)}}\otimes\cdots\otimes v_{i_{\sigma(m)}}.
\end{eqnarray*}
Notice that this action of the permutation group commutes with the
induced action of $G=SU(n)$ on $T^m(V)$, and therefore we have a
representation $\tau_m$ of $G\times {\cal S}_m$ on $T^m(V)$. It
also follows that image of $T^m(V)$ under a Young symmetrizer is
invariant under $G$.  It is well-known that

\begin{thm}
\label{thm:Young2} Every partition $T:m=m_1+\cdots +m_k$ with
$m_1\ge m_2\ge\cdots\ge m_k$ determines a unique irreducible
representation $\lambda_T$ of ${\cal S}_n$, and every irreducible
representation of ${\cal S}_n$ is of the form $\lambda_T$. The
degree of $\lambda_T$ is the number of standard Young tableaux
whose underlying Young diagram is $T$.
\end{thm}

The basic result of Schur-Weyl theory can be summarized as
follows:

\begin{thm}
\label{thm:Young1} The representation $\rho_T$ of $G$ is
irreducible. For every Young diagram $T$ corresponding to a
partition of $m$, let $Z_T\subset T^m(V)$ be the minimal linear
subspace containing ${\rm Im}{\sf{C}}_T$ and invariant under
action of $SU(n)\times {\cal S}_m$. $Z_T$ has dimension $\deg
(\rho_T)\deg (\lambda_T)$ and is irreducible under the
representation $\tau_T=\rho_T\otimes\lambda_T$ of $SU(n)\times
{\cal S}_m$. $\deg (\rho_T)$ is equal to the number of
semi-standard Young tableaux whose underlying diagram is $T$.
Furthermore $T^m(V)$ admits of the decomposition,  as a $G\times
{\cal S}_m$-module (under $\tau_m$),
\begin{eqnarray*}
T^n(V)\simeq \sum_T Z_T,
\end{eqnarray*}
where the summation is over all partitions of $T$ of $m$ with
$k\le n$ parts.  Let ${\cal A}_T(G)$ and ${\cal A}_T({\cal S}_m)$
denote the algebras of linear transformations of $Z_T$ generated
by the matrices $\rho_T(g)\otimes I$, ($g\in G$), and $I\otimes
\lambda_T(\sigma)$, ($\sigma\in {\cal S}_m$). Then the full matrix
algebra on $Z_T$ has the decomposition ${\cal A}_{T}(G)\otimes
{\cal A}_{T}({\cal S}_m)$.
\end{thm}

An irreducible representation $\rho$ of $SU(n)$ occurs in
$L^2(G/K)$ if and only if it has a $K$-fixed vector (Frobenius
reciprocity), and since $(G,K)$ is a symmetric pair the space of
$K$-fixed vectors is one dimensional.  We have

\begin{prop}
\label{prop:jshaffaf1} Let $G=SU(n)$ and $K=S(U(k) \times U(n-k))$
where $k\ge n-k$. An irreducible representation $\rho$ of $G$ has
a $K$-fixed vector if and only if the corresponding Young diagram
is of the form

 \baselineskip=.04in
$
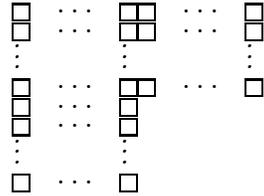
\begin{array}{llllll}
~~~~~~~~~~~~~~~~~~~~~~~~~~~~~~~~~~~~~~~~~
&\fbox{}&\ldots&\fbox{}\fbox{}&\ldots&\fbox{}\\[-0.25cm]
&\fbox{}&\ldots&\fbox{}\fbox{}&\ldots&\fbox{}\\[-0.2cm]
&\vdots&&\vdots&&\vdots\\[-0.25cm]
&\fbox{}&\ldots&\fbox{}\fbox{}&\ldots&\fbox{}\\[-0.25cm]
&\fbox{}&\ldots&\fbox{}&&\\[-0.25cm]
&\fbox{}&\ldots&\fbox{}&&\\[-0.2cm]
&\vdots&&\vdots&&\\[-0.25cm]
&\fbox{}&\ldots&\fbox{}&&\\
\end{array}
$
\vskip .3 cm \noindent where there are $k$ squares in the first
column and $n-k$ squares in the last, and the number of columns of
lengths $k$ and $n-k$ are equal.
\end{prop}

\noindent {\bf Proof} - Let $T$ be a Young diagram of the form
specified in the lemma and let $r$ denote the number of columns of
length $k$ (or $n-k$) , and $\{T\}$ denote the semi-standard Young
tableau where the first $r$ columns are filled with integers
$1,\ldots,k$ and the last $r$ columns are filled with integers
$k+1,\ldots,n$. The action $g\in U(k)\times U(n-k)$ on the vector
$v_T$ corresponding to $\{T\}$ is given by
\begin{eqnarray*}
v_T\longrightarrow (\det g_1 \det g_2)^r,~~~~ where~~
g=\begin{pmatrix} g_1&0\cr o&g_2
\end{pmatrix}
\end{eqnarray*}
Therefore $v_T$ is fixed by $K$.  The converse statement that the
existence of $K$-fixed vector implies the corresponding Young
diagram is of the required form will be proven only for $k=n-1$
which is the only case needed here.  We make use of the following
simple

\begin{lemma}
\label{lem:jshaffaf2} Let ${\mathcal T}\subset SU(n)=G$ be the
maximal torus of diagonal matrices.  An irreducible representation
$\rho_T$ of $G$ contains a $T$-fixed vector if and only if the
corresponding Young diagram contains $rn$ squares for some
positive integer $r$, and the fixed vector is represented by a
Young tableau with the same number of 1's, 2's,$\ldots, n$'s.
\end{lemma}

\noindent {\bf Proof of Lemma} - Let $e^{i\tau_j}$ denote the
diagonal entries of a matrix in ${\mathcal T}$.  Then the only
relation among $\tau_j$'s is $\sum\tau_j=0$.  Therefore there is a
vector in the representation space of $\rho_T$ fixed by
$\mathcal{T}$ if and only if there is a semi-standard Young
tableau $\{T\}$ containing the same number of 1's, 2's,$\ldots,
n$'s which proves the lemma.

Let $\rho_T$ denote the irreducible representation $U(n)$ and the
Young diagram $T$ correspond to the partition $m=m_1+\ldots +m_n$.
According to the Branching Law [Kn, page 569] the restriction of
$\rho_T$ to $U(n-1)$ decomposes into a direct sum of irreducible
representations with multiplicity one according to partitions
$l=l_1+\ldots +l_{n-1}$ such that
\begin{equation}
\label{eq:jshaffaf11} m_1\ge l_1\ge m_2\ge l_2\ge \ldots \ge
l_{n-1}\ge m_n.
\end{equation}
In order for the restriction of $\rho_T$ to $U(n-1)$ to contain
the representation $\det^r$ it is necessary and sufficient that
one of the partitions of $l$ be of the form
\begin{eqnarray*}
l=l^\prime+l^\prime+\ldots +l^\prime, ~~~{\rm
that~is}~~l=rl^\prime.
\end{eqnarray*}
Therefore by (\ref{eq:jshaffaf11}), the representation $\det^r$
occurs in the restriction of $\rho_T$ to $U(n-1)$ if and only if
\begin{equation}
\label{eq:jshaffaf12} l^\prime = m_2=m_3=\ldots =m_{n-1}.
\end{equation}
Since for irreducible representation of $SU(n)$ it is only
necessary to consider Young diagrams with $n-1$ rows, it follows
from (\ref{eq:jshaffaf12}) that the restriction of an irreducible
representation $\rho_T$ of $SU(n)$ to $U(n-1)$ contains $\det^r$
if and only if
\begin{eqnarray*}
m_1\ge m_2=m_3=\ldots =m_{n-1}.
\end{eqnarray*}
By Lemma \ref{lem:jshaffaf2} such a representation contains a
${\mathcal T}$-fixed vector if and only if there is a Young
tableau $\{T\}$ with the same number of 1's,2's,$\ldots, n$'s.
Furthermore the one dimensional invariant subspace transforming
according to $\det^r$ under $U(n-1)$ is spanned by the Young
tableau $\{T\}$ where the number of $1$'s and $n$'s in the first
row is $l^\prime =m_2$.  Therefore $m_1=2l^\prime$ and the proof
of the Proposition is complete.  $\blacksquare$ \vskip .3 cm
\noindent {\bf Proof of Theorem \ref {thm:jshaffaf2} }- In order
to prove the theorem we recall the relevant aspect of the
Plancherel theorem for a compact connected semi-simple Lie group.
The Fourier transform of the spherical measures $\lambda_a$ is
\begin{equation}
\label {eq:shaffaf8} F(\lambda_a)(\rho)
=\widehat{\lambda_a}(\rho)= \int_{O_a} \rho(x) d\lambda_a =
\phi_{\rho}(a)\textrm{vol}(O_a) ~,
\end{equation}
with similar expression for $\widehat{\lambda_b}$, where
$\phi_{\rho}$  is the elementary $K$-spherical function
corresponding to the irreducible representation $\rho$ containing
a $K$-fixed vector. It is well-known that every elementary
$K$-spherical function on $G$ is of the form $\rho_{11}(g)$ in
which $\rho_{11}$ is the $(11)$-entry of the matrix of $\rho$
relative to an orthonormal basis $v_1,\ldots,v_N$ where $\rho
(K)v_1=v_1$ (see [H2], page 414).  According to Proposition
\ref{prop:jshaffaf1} irreducible representations of $SU(n)$
containing a $K$ fixed vector are parameterized by integers $m$
corresponding to partitions
\begin{eqnarray*}
N=2m~+~m~+~m~+~\ldots~+~m~~~~~(n-1)~{\rm summands}.
\end{eqnarray*}
We need

\begin{lemma}
\label{lem:shaffaf61} Let $\rho$ be a spherical representation of
the group $G$ (i.e., containing a $K$-fixed vector $v$), then the
corresponding elementary spherical function
\begin{eqnarray*}
|\phi_m (\exp tH)|=
C~\frac{(n-1)^{n-\frac{1}{2}}m^{m+\frac{1}{2}}}{(n+m-1)^{n+m-\frac{1}{2}}}~\frac{t^{-n+\frac{1}{2}}}{\sqrt{m}}.
\end{eqnarray*}
for some constant $C$ independent of $m$ and $n$.
\end{lemma}

\noindent {\bf Proof} - This lemma is probably well-known to
experts in spherical functions but since the author does not know
of a specific reference for it in the form suitable for this work,
a proof is sketched. However, in [SC] estimates for Jacobi
polynomials are used to establish precise rates of convergence for
certain diffusion processes.  By applying the radial part of the
Laplacian to the elementary spherical function $\phi_m$, one
obtains a second order linear ordinary differential equation with
regular singular points for it. Consequently one obtains [H2]
\begin{eqnarray*}
\phi_m (\exp tH)= F(m+n,-m,n;\sin^2t)
\end{eqnarray*}
where $F$ is the hypergeometric function which reduces to the
Jacobi polynomial $P^{n-1,n}_m (\frac{\cos^2t}{2})$ (except for
normalization by a constant) (see [AAR] for explanation of
notation and extensive treatment of Jacobi polynomials).  Now
\begin{eqnarray*}
P^{n-1,n}_m(1)=\frac{(n+m-1)!}{(n-1)!m!}.
\end{eqnarray*}
Since $\phi_m(e)=1$, $P^{n-1,n}_m$ should be normalized
accordingly.  Estimates for Jacobi polynomials are obtained by
examining the behavior of their generating function on the unit
circle and applying standard methods for obtaining estimates from
generating functions. In fact one obtains (see [AAR] especially
page 350)
\begin{equation}
\label{eq:shaffaf61} P^{n-1,n}_m(\cos\theta)=\begin{cases}
\theta^{-n+\frac{1}{2}} O(\frac{1}{\sqrt{m}}),~~{\rm
for}~~\frac{c}{n}\le \theta\le \frac{\pi}{2};\\
O(m^{n-1}),~~{\rm for}~~0\le \theta\le\frac{c}{n};
\end{cases}
\end{equation}
for a suitable constant $c$ as $m\to\infty$.  Substituting in the
expression for $\phi_m$ in terms of Jacobi polynomials we obtain
the desired estimate.  $\blacksquare$

\noindent Because  $\phi_m$ is a spherical function, the value of
this function on the spherical class $O_a$ is constant and so we
can take $a=\exp(t_1H)$.  Now by (\ref{eq:shaffaf8}) for the
Fourier transform of the spherical measure $\lambda_a$ and
$\lambda_b$ we obtain:

\begin{eqnarray*}
F(\lambda_a)(\rho_m) = \phi_m(a) \textrm{vol}(O_a) = \phi_m(
\exp t_1H) \textrm{vol}(O_a)\\
F(\lambda_b)(\rho_m) = \phi_m(b) \textrm{vol}(O_b) = \phi_m( \exp
t_2H) \textrm{vol}(O_b) ~,
\end{eqnarray*}

\noindent where $\rho_m$ is the spherical representation
corresponding to $m$.

Now applying the Plancherel theorem for the function $(\lambda_a
\star \lambda_b)^{l(n)}$ we obtain
\begin{eqnarray*}
 \| (\lambda_a \star \lambda_b)^{l(n)} - 1 \|^2_{L^2}
&=& \ \sum_{m>0} d_{\rho_m} (\phi_m(t_1) \phi_m(t_2))^{2l(n)} ~,
\end{eqnarray*}
We want to find $l(n)$ such that the sequence
 \begin{equation}
 \label{eq:shaffaf11}
 c_n=\sum_m
d_{\rho_m} (\phi_m(t_1) \phi_m(t_2))^{2l(n)}
\end{equation}
converges to zero when $n$ goes to infinity.  This will ensure
that $(\lambda_a \star \lambda_b)^{l(n)}$ converges to the Haar
measure in $L^p$-norm for $1\le p\le 2$ as $n\to \infty$ since on
compact groups of fixed finite volume $L^2$ norm dominates $L^p$
for $p\le 2$.

For analyzing the sequence $\{c_n \}$ we need to compute the
dimension of the representation $\rho_m$. By Weyl's dimension
formula the dimension of the irreducible representation $\rho_T$
of $U(n)$ determined by the Young diagram $ T: m = m_1 + m_2 +
\ldots +m_n$ is $$ \frac{\mathcal{D}(a_1,a_2, \ldots ,
a_n)}{\mathcal{D}(m-1,m-2, \ldots , 0)}~ ~ ,$$ where $a_k = m_k +
n - k$ and $\mathcal{D}(a_1,a_2, \ldots , a_n) = \prod_{j<k} (a_j
- a_k)$. In our situation the spherical representation $\rho_m$
has the Young diagram characterized  in Proposition
\ref{prop:jshaffaf1} and the corresponding partition is $ m n = 2m
+ m + \ldots + m $ where $2m$ is the number of columns in the
corresponding Young diagram.  Therefore
\begin{eqnarray*}
a_1 = 2m+n-1,~~ a_2 = m+n-2,~ a_3 = m+n-3, \ldots ,a_{n-1} = m+1,
~ a_n=0
\end{eqnarray*}
By the Weyl's dimension formula the dimension $d(m,n)$ of the
representation $\rho_m$ is
\begin{eqnarray*}
d(m,n) &=& (2m+n-1)\frac{\prod_{k=2}^{n-1} [(n-k-1)!
(m+n-k)^2]}{(n-2)!(n-1)! \ldots 2! 1!}\\
&=& (2m+n-1)\frac{\prod_{k=2}^{n-1} (m+n-k)^2}{(n-2)!(n-1)!} \\
&=& \frac{((m+n-2)!)^2}{(m!)^2 (n-2)! (n-1)!} (2m +n-1).
\end{eqnarray*}
Applying the Sterling estimate $ n!\sim \frac{1}{\sqrt{2 \pi}}
n^{n+\frac{1}{2}} e^{-n}$ we obtain
\begin{eqnarray*}
d(m,n) \sim \frac{e}{2 \pi} \frac{(m+n-2)^{2(m+n)-3}}{m^{2m+1}
n^{2n-2}}(2m+n-1)
\end{eqnarray*}
Substituting for $\phi_m(t)$ from the Lemma \ref{lem:shaffaf61} we
obtain
\begin{eqnarray*}
c_n = b^{(2n-1)l(n)} \sum_{m=1}
\frac{(m+n-2)^{2(m+n)-3}(2m+n-1)}{m^{2m+1} n^{2n-1}} \bigg(
\frac{(n-1)^{2n-1} m^{2m}}{(n+m-1)^{2(m+n)-1}} \bigg)^{2l(n)}
\end{eqnarray*}
where $b = \frac{1}{t_1 t_2}$.  Now we decompose the summation
$c_n$ into two parts $s_1$ and $s_2$ as follows:

\begin{eqnarray*}
c_n = b^{(2n-1)l(n)}(s_1 + s_2) = b^{(2n-1)l(n)}(\sum_{m\le n-1} ~
+ ~ \sum_{m > n-1})
\end{eqnarray*}
 Since in $s_1$ the summation is over $m \leq n-1$ we have
 \begin{eqnarray*}
 s_1 &=& \sum_{m \leq n-1}
 (m+n-1)^{2(m+n)-3-2l(n)(2(m+n)-1)}(2m+n-1) m^{4m l(n)-2m-1} n^{2(2n-1)l(n)-2n+2}  \\
 &\leq & n^{2(2n-1)l(n)-2n+2} \sum _{m \le n-1} (2(n-1))^{2(m+n)-3-2l(n)(2(m+n)-1)}
 (n-1)^{4ml(n)-2m-1}
  \\ &= & 2^{2n-4nl(n)+2l(n)-3}(n-1)^{2n-4nl(n)+2l(n)-3}
  n^{2(2n-1)l(n)-2n+2}.\\&&
 \sum_{m \leq n-1} 2^{2m-4ml(n)} (n-1)^{2m-4ml(n)} (n-1)^{4ml(n)-2m}
\\ & = & 2^{2n-4nl(n)+2l(n)-3} n^{-1}  \sum_{m \leq n-1}
 (4^{1-2l(n)})^m
 \end{eqnarray*}
 Therefore
 \begin{eqnarray*}
 b^{(2n-1)l(n)} s_1 \le 2^{2n-4nl(n)+2l(n)-3} n^{-1} b^{(2n-1)l(n)}
 \frac{a^n-1}{a-1},
\end{eqnarray*}
where $a = 4^{1-2l(n)}$.  Now for $l(n)\ge C_1\log n$ and $C_1$
sufficiently large and depending on $t_1$ and $t_2$, we obtain
\begin{equation}
\label{eq:jshaffaf1} s_1\le \frac{C_2}{n^\epsilon},
\end{equation}
for some $\epsilon>0$ and some constant $C_2$ depending only on
$t_1$ and $t_2$.

Now we estimate $s_2$ where the summation is over $m > n-1$.  It
is clear that

\begin{eqnarray*}
s_2 &\leq &\sum_{m > n-1} 3m (2m)^{2(m+n)-3-2l(n)(2(m+n)-1)}
m^{4ml(n)-2m-1} n^{2(2n-1)l(n)-2n+2}
 \\ &\ = & b^{(2n-1)l(n)} 2^{2n-4nl(n)+2l(n)-3} n^{2(2n-1)l(n)-2n+2} \sum_{m>n-1} m^{2n +2l(n)-4n l(n)-3}2^{2m-4ml(n)}
 \\ & \le & b^{(2n-1)l(n)} 2^{2n-4nl(n)+2l(n)-3} n^2
 \sum_{m=n}^{\infty} a^m \\ & = & b^{(2n-1)l(n)} 2^{2n-4nl(n)+2l(n)-3} n^2
 \frac{a^n}{1-a},
\end{eqnarray*}
where $a = 4^{1-2l(n)}$ as before.  Consequently for $l(n)\ge C_3
\log n$ and $C_3$ sufficiently large and depending on $t_1$ and
$t_2$ we have
\begin{equation}
\label{eq:jshaffaf2} s_2\to 0~~~{\rm as} ~ n\to\infty.
\end{equation}
Therefore $c_n = s_1+s_2$ tends to zero and the proof of the
theorem is complete. $\blacksquare$

\vfil\eject

\vskip .9 cm

\noindent Institute for Studies in Theoretical Physics and
Mathematics, Tehran, Iran, {\it and} \\ Sharif University of
Technology, Tehran, Iran. \vskip .4 cm

\end{document}